\documentclass{amsart}[11pt,leqno] \usepackage{amssymb,enumerate}

\numberwithin{equation}{section}

\theoremstyle{plain}
\newtheorem{theorem}[equation]{Theorem}
\newtheorem{lemma}[equation]{Lemma}
\newtheorem{corollary}[equation]{Corollary}
\newtheorem{proposition}[equation]{Proposition}
\theoremstyle{definition}
\newtheorem{Def}[equation]{Definition}
\newtheorem{Ex}[equation]{Example}
\newtheorem{Rem}[equation]{Remark} 
\newtheorem{Not}[equation]{Notation}

\newcommand{\F}{{\mathbb F}}
\newcommand{\Q}{{\mathbb Q}}

\newcommand{\Z}{{\mathbb Z}}
\newcommand{\Mu}{{\boldsymbol \mu}}

\newcommand{\ga}{{\mathfrak a}}
\newcommand{\gD}{{\mathfrak D}}
\newcommand{\gf}{{\mathfrak f}}
\newcommand{\gl}{{\mathfrak l}}
\newcommand{\gn}{{\mathfrak n}}
\newcommand{\go}{{\mathfrak o}}
\newcommand{\gp}{{\mathfrak p}}
\newcommand{\gS}{{\mathfrak S}}

\newcommand{\cA}{{\mathcal A}}
\newcommand{\cD}{{\mathcal D}}
\newcommand{\cI}{{\mathcal I}}
\newcommand{\cO}{{\mathcal O}}
\newcommand{\cS}{{\mathcal S}}
\newcommand{\cV}{{\mathcal V}}
\newcommand{\cW}{{\mathcal W}}

\newcommand{\tE}{\tilde{E}} 
\newcommand{\tF}{\tilde{F}}

\newcommand{\tL}{{\tilde{L}}}
\newcommand{\tM}{{\tilde{M}}}
\newcommand{\Gal}{{\rm Gal}}

\newcommand{\Image}{{\rm Image} \, }
\newcommand{\GL}{{\rm GL}}
\newcommand{\SL}{{\rm SL}}
\newcommand{\GS}{{\rm GSp}}
\newcommand{\SP}{{\rm Sp}}
\newcommand{\ord}{{\rm ord}}
\newcommand{\rank}{{\rm rank \, }}
\newcommand{\Aut}{{\rm Aut}}
\newcommand{\End}{{\rm End}}
\newcommand{\Hom}{{\rm Hom}}

\begin{document}

\title{Arithmetic of division fields}

\author[A. Brumer]{Armand Brumer}
\address{Department of Mathematics, Fordham University, Bronx, NY 10458}
\email{brumer@fordham.edu}
\author[K. Kramer]{Kenneth Kramer}
\address{Department of Mathematics, Queens College and the Graduate Center (CUNY), 65-30 Kissena Boulevard, Flushing, NY 11367, Flushing, NY 11367}
\email{kkramer@gc.cuny.edu}

\thanks{Research of the second author partially supported by NSF grant DMS 0739346}


\date{February 21, 2011}

\begin{abstract}
We study  the arithmetic of division fields of semistable abelian varieties $A_{/\Q}.$  The Galois group of $\Q(A[2])/\Q$ is analyzed when the conductor is odd and squarefree.    The irreducible semistable mod 2 representations of small conductor are determined under GRH.  These results are used in {\em Paramodular abelian varieties of odd conductor}, arXiv:1004.4699.
\end{abstract}


\maketitle

\section{Introduction}   
This note contains results needed in \cite{BK2} and of independent interest.  We write $S$ for a set of primes, $N_S$ for their product and $\ell$ for a prime not in $S$.  If $F/\Q$ is Galois, $\cI_v(F/\Q)$ denotes the inertia group at a place $v$ of $F$.

\begin{Def}[\cite{BK1}] \label{con}  The Galois extension $F/\Q$ is $(\ell,N_S)$-{\em controlled} if
\begin{enumerate}[{\rm i)}]
\item $F/\Q$ is unramified outside $S \cup \{\ell, \infty\}$;
\item $\cI_v(F/\Q) = \langle \sigma_v \rangle$ is cyclic of order $\ell$ for all ramified $v$ not over $\ell$;
\item $\cI_\lambda(F/\Q)^u = 1$ for all $u > 1/(\ell-1)$ and $\lambda$ over $\ell$, using the upper numbering of Serre as in  \S \ref{FontSect}.
\end{enumerate}
\end{Def}

We denote by $V$  a finite dimensional vector space over the finite field $\F$ of characteristic $\ell$ with  $q = |\F|$. Additional structure on $V$, such as a symplectic pairing or Galois action, is often imposed. 
  
\begin{Def} \label{SSRep}  Let $V$ be an $\F[G_\Q]$-module and $F = \Q(V)$.  The set $S$ of rational primes $p \ne \ell$ ramified in $F/\Q$ comprises the {\em bad primes} of $V$.  Declare $V$   {\em semistable} if $F$ is $(\ell,N_S)$-controlled and $(\sigma_v-1)^2(V) = 0$ for all $v$ lying over the primes of $S$.
\end{Def}

Throughout, $A_{/\Q}$ is a semistable abelian variety with good reduction at $\ell$ and $\End_\Q A=\go$ is the ring of integers in a totally real number field.  If $\gl$ is a prime over $\ell$ in $\go$ and $\go/\gl = \F$, then $V=A[\gl]$ is semistable \cite{Gro, Fo}.   The conductor of $A$ has the form $N_A = N^d$ with $d = [\go\!:\!\Z]$.  Since inertia over each bad prime $p$ is tame, 
\begin{equation} \label{condexp}
\ord_p(N) = \dim_\F V/V^\cI = \dim_\F \, (\sigma_v-1) V.
\end{equation}
In \S\ref{grpthy}, we use known results on symplectic representations generated by transvections to describe $\Gal(\Q(W)/\Q)$  for  constituents $W$ of $V$ with squarefree conductor, assuming $\gl$ lies over 2.

A {\em stem field} for a Galois extension $F/k$ is an intermediate field $K$ whose Galois closure over $k$ is $F$.  If $G = \Gal(F/k)$ acts faithfully and transitively on a set $X$, the fixed field of the stabilizer $G_x$ of any $x$ in $X$ is a stem field.   A formula for the discriminant $d_{K/k}$ is given  in \S \ref{DiscSect}  and applied to semistable Galois modules.   By relating number-theoretic properties of $K$ and $F$, certain computations become feasible, since $K$ has smaller degree and discriminant than $F.$  

Suppose $E/\Q_\ell$ is a Galois extension of $\ell$-adic fields satisfying (\ref{con}iii).  In \S \ref{FontSect}, we find conditions on the ray class conductor of an abelian extension $L/E$ so that (\ref{con}iii) also holds for the Galois closure of $L/\Q_\ell.$  
The maximal $(2,N)$-controlled extension for all odd $N \! \le \! 79$ and for $N\!=\!97,$ is determined in \S \ref{Odl}, thanks to  \S \ref{FontSect} and Odlyzko's GRH bounds. We also construct a $(2,127)$-controlled extension of degree 161280 with root discriminant just above the asymptotic Odlyzko bound, but finiteness of a maximal one  is unknown.

A finite flat group scheme $\cV$  over $\Z_\ell$ admits a filtration  $0 \subseteq \cV^m \subseteq \cV ^0 \subseteq \cV$ with connected component $\cV ^0$, \'etale quotient $\cV ^{et} = \cV /\cV ^0$, multiplicative subscheme $\cV ^m$ and biconnected subquotient $\cV ^b = \cV ^0/\cV ^m$.   Let $\lambda$ be a place over $\ell$ in $F=\Q(V)$  and  $\cD_\lambda$ its decomposition group. We denote the corresponding $\F[\cD_\lambda]$-modules of $F_\lambda$-valued points by $V $, $V ^{et}$, $V ^m$ and $V ^b$, respectively.

\begin{Def}[\cite{BK2}]
$A_{/\Q}$ is $\go$-{\em paramodular} if $\dim A=2d$, with $d = [\go\!:\!\Z]$.
\end{Def}

Let $A$ be $\go$-paramodular, with $\go/\gl \simeq \F_2.$  When $A[\gl]$ is irreducible, estimates for the discriminant of a stem field of $\Q(A[\gl])$ are obtained in \S \ref{ParaSect}.  The reducible case leads to ray class fields whose conductors are controlled by the results of \S \ref{FontSect}.  This information depends on the structure of $A[\gl]$ as a group scheme and is used in \cite{BK2}. 

\medskip

Concluding questions and comments appear in \S \ref{spec}.

\section{Mod 2 representations generated by transvections} \label{grpthy}
A {\em transvection} on $V$ is an automorphism of the form $\tau(x)=x+\psi(x) \, z$, with  $\psi\!: V \to \F$ a nonzero linear form and $z\ne 0$ in $\ker \psi.$  Assume $V$ admits a non-degenerate alternating pairing $[ \, , \, ]: \, V \times V \to \F$  preserved  by  $\tau$ and let $\dim V = 2n$. Then $\tau(x)=x + a \, [z,x] z$ for some $z\in V$ and $a\in\F^\times$.  When $a$ is a square in $\F$, we may take $a = 1$.  For $x$ and $z$ in $V$, define $\tau_{[z]}$ by
\begin{equation}\label{sx}
\tau_{[z]}(x) = x + [z,x] z .
\end{equation}
Assume  that  $\ell = 2$ for the rest of this section, unless otherwise noted.   

A quadratic form $\theta$ on the symplectic space $V$ is called a {\it theta characteristic} if 
$\theta(x+y) = \theta(x)+\theta(y) + [x,y]$ for all $x,y$ in $V$.  Theta characteristics form a principal homogeneous space over $V$, with 
$
(\theta+a)(x) = \theta(x) + [a,x]^2
$
for $a$ in $V$.  We identify $a$ with $[\, a, - ]$ under the Galois isomorphism $V \simeq \Hom_\F(V,\F)$.  Elements  $\sigma$ in $\SP(V)$ act  by $\sigma(\theta)(x) = \theta(\sigma(x)).$ Then $\sigma(\theta+a) = \sigma(\theta) + \sigma(a)$ and
\begin{equation} \label{tQ}
\tau_{[z]}(\theta) = \theta + \sqrt{1+\theta(z)} \, z.
\end{equation}

Fix a  symplectic basis $\{e_1, \dots, e_{2n}\}$ for $V$ with $[e_i, e_j]=1$ if $|i-j|=n$ and 0 otherwise.  Let  $\wp(x)=x^2-x$ be the Artin-Schreier function. 
Depending on whether or not the Arf invariant ${\rm Arf}(\theta)=\sum_i \theta(e_i) \theta(e_{i+n})$ vanishes in $\F/\wp(\F),$ we say $\theta$ is even or odd and write $O_{2n}^{\pm}$  for the corresponding orthogonal group.
Further, $\SP(V)$ acts transitively on  the sets $\Theta_{2n}^\pm$ of even and odd characteristics and 
\begin{equation}\label{thetacount}
|\Theta_{2n}^\pm|=  \frac{1}{2}q^n(q^n \pm 1).
\end{equation}

Denote the symmetric, alternating, dihedral and cyclic groups by $\cS_n,$ $\cA_n$, $D_n$, $C_n$ respectively.

\begin{proposition}[\cite{McL}] \label{McL} 
If $\F=\F_2$ and  $G\subsetneq\SL(V)$ is an   irreducible subgroup  generated by transvections, then $\dim V =2n$ with $n\ge 2$ and $G$ is  $O^{\pm}_{2n}(\F_2)$, $\SP_{2n}(\F_2)$ or  $\cS_m$ with  $2n+1\le m \le 2n+2$.
Also, $G$  has trivial center and is self-normalizing in $\SL(V)$.
\end{proposition}

\begin{proposition} \label{trgp} 
Let $V$ be a symplectic space of dimension $2n$. An irreducible subgroup $G$ of $\SP(V)$ generated by transvections is  one of the  
following, with $\F'\subseteq\F$:
\begin{enumerate}[{\rm i)}]
\item   dihedral, $D_m$  with $m$ dividing one of $|\F| \pm 1 $ and $n=1$;
\item  orthogonal, $O^{\pm}_{2n}(\F')$ for  $n\ge 2$;
\item  symplectic, $\SP_{2n}(\F')$; 
\item  symmetric, $\cS_m$ for $n\ge 2$ and $2n+1\le m \le 2n+2$.
\end{enumerate}
Moreover, $G$  has trivial center and is self-normalizing in $\SP(V)$.
\end{proposition}

\proof 
If $V$ is imprimitive, then $V$ is monomial \cite{Zal1}, say $V = {\rm Ind}_{H_1}^G(V_1)$, with $V_1=\F e_1$ and $[G:H_1]=\dim V = 2n$.  Arrange that $G= \cup \, g_iH_1$, with $g_1=1$ and $V_i=g_i(V_1)=\F e_i$, and let $\pi: G \to \cS_{2n}$ by $g V_i=V_{\pi(g)i}$.  Since $\pi(G)$ is transitive  and generated by transpositions, namely the images of the transvections, $\pi(G)=\cS_{2n}$.  For $h$ in $H = \ker \pi$, we have $he_i=\chi_i(h)e_i$ and so the pairing on $V$ satisfies $[e_i,e_j] = [he_i,he_j] = \chi_i(h)\chi_j(h)[e_i,e_j]$.  Hence $[e_i,e_j]=0$ or $\chi_i(h)\chi_j(h)=1$.  Because the pairing is perfect and $\pi(G)$ is doubly transitive, we must have $[e_i,e_j]  \ne 0$ and $\chi_i(h)\chi_j(h)=1$ for all $i\neq j$.  If $n \ge 2$, then $\chi_i(H) = 1$ for all $i$,  $H = 1$ and $\pi$ is an isomorphism.  The stabilizer $H_1$ of $V_1$ is isomorphic to  $\cS_ {2n-1}$  and so the character $\chi_1: H_1 \to \F^\times$ is trivial.  Since $\sum \, g_i(e_1)$ is a non-trivial fixed point, $V $ is reducible.  Now combine  \cite[Ch. II, \S 8.27]{HBI} and  \cite{Kan,KM} to get our list.  

If  $g$ in $\SP_{2n}(\F)$ normalizes $G$ and $\sigma$ is in $\Gal(\F/\F')$, then $g^{\sigma}g^{-1}$ centralizes $G.$  Our representations  are absolutely irreducible and the center of $\SP_{2n}(\F)$ is trivial, so $g$ is in $\SP_{2n}(\F')$.   To verify that the center is trivial and $G = \cS_m$ is self-normalizing in $\SP_{2n}(\F_2)$ when $m \ne 6$,  use the fact that all automorphisms are inner and absolute irreducibility.  Note that $\cS_6 \simeq \SP_4(\F_2)$.  The dihedral case is easily checked.  See \cite{Dye} for the other  
cases.   \qed

\begin{Rem} \label{SnRep} 
As to (iv) above, note that  $\cS_m$ acts by permutation on $$Y=\{(a_1,..,a_m) \in \, \F_2^m\ |\  a_1 + \dots + a_m =0\}$$ with pairing $[(a_i),(b_i)]=\sum a_i  b_i$.  Let $V=Y/\langle (1,\dots,1)\rangle$ or $V=Y$ according as $m$ is even or odd.   Then $V$ is  irreducible and transpositions in $\cS_m$ correspond to transvections on $V$.  This action of $\cS_m$  and that of Galois on $J[2]$ for a hyperelliptic Jacobian are compatible.
\end{Rem}

\begin{lemma} \label{Cliff}
Let $V$ be an irreducible $\F[G]$-module and let $P$ be the subgroup of $G$ generated by transvections. If  $P$ is not trivial, then $V_{|P}$ is the  direct sum of $r$ irreducible $\F[P]$-modules  $W_i$ and $P = Q_1 \cdots  Q_r$ is a direct product, with
$
Q_i = \langle \sigma \in P \, | \, \sigma_{|W_i}  {\rm \; is \; a \; transvection \; and \;} \sigma_{|W_j} = 1 {\rm \; for \, all \;} j \ne i \rangle.
$
If $V$ is symplectic, then the $W_i$ are symplectic and the sum is orthogonal.
\end{lemma}

\proof
Since  $P$  is  normal, Clifford's theorem applies. Let $W_1$ be an irreducible submodule of  
$V_{|P}$, 
$
H=\{ h {\rm \,  in \, } G \, | \, h(W_1) \simeq W_1 {\rm \; as \; } P{\mbox-}{\rm module} \} 
$ and $X = \sum_{h \in   H} h(W_1)$.  Then $V = {\rm ind}_H^G(X)$  and $X_{|P} \simeq eW_1$ is isotypic. If $G = \cup_1^r \, g_iH$ is  a coset decomposition with  
$g_1=1,$  then $V_{|P} \simeq \oplus_1^r \,  eW_i$ with $W_i = g_i(W_1)$. 
For any  transvection $\tau$, we have
$
1 = \dim \, (\tau-1)(V) = e \, \sum_1^r \dim \, (\tau-1)(W_i).
$
Thus $e=1$ and  $\tau$ is in $Q_i$ for a unique  
index $i$.  
Moreover $Q_i = g_i Q_1 g_i^{-1}$ is normal in $P$ and  $P = Q_1  
\cdots Q_r$ is a direct product.  

Now suppose $V$ symplectic and $\tau$ a transvection in $Q_i$.  Then $(\tau-1) W_i=\langle z\rangle$ with $z$ in $W_i \cap W_j^\perp$ for all $j \ne i$, but  not in $W_i^\perp$.  Irreducibility of $W_i$ implies that  $W_i\subseteq W_j^\perp$ and $W_i \cap W_i^\perp = 0$.  Hence $W_i$ is symplectic. \qed

\medskip
\begin{proposition}\label{exc2} 
Let $V$ be an irreducible symplectic $\F[G_\Q]$-module with squarefree conductor $N$ and let $F=\Q(V)$.  Let $P$ be the subgroup of $G = \Gal(F/\Q)$ generated by transvections. If $P=G$, then  $G$ is as in {\rm Prop.\! \ref{trgp}}. 

Otherwise, $V={\rm ind}_P^G W$ and $G \simeq Q \wr C_2$, where $Q$ is in the list in {\rm Prop.\! \ref{trgp}}.  Moreover $F^P = \Q(i)$ and  $N=\gn \overline{\gn}$ in $\Z[i]$, where $\gn$ generates the conductor ideal of $W$ as $\F[G_{\Q(i)}]$-module.
\end{proposition}

\proof  
Since $\ord_{p_v}(N) = 1$, any generator $\sigma_v$ of  $\cI_v(F/\Q)$ is a transvection.  Prop.\! \ref{Onlyell} shows that the fixed field $F^P = \Q(i)$.  The restriction $V_{|P}$ is  reducible by Lemma \ref{trgp} and so  $V$ is induced. Hence $H=P \simeq Q_1 \times Q_2$ and $G \simeq Q_1 \wr C_2$  is a wreath product, thanks to Lemma \ref{Cliff}.  The conductor formula for an induced module gives $N = \gn \overline{\gn},$ where $\gn \in \Z[i]$ is the  odd part of the Artin conductor of $W$, since $\Q(i)$ is unramified at odd places.  \qed

\medskip

\begin{Rem}
In  Prop.\! \ref{exc2}, if we take $\F=\F_2$ but do not assume $V$ symplectic,   the conclusions obtain, with ``{\rm Prop.\! \ref{trgp}}'' replaced by ``{\rm Prop.\! \ref{McL}}." 
\end{Rem}

\begin{Rem} \label{tdelta}
The conjugacy class of any involution $\sigma$ in $\SP(V)$ has invariants  $t = \rank  (\sigma-1)$ and $\delta$, with $\delta=0$ if $[v,(\sigma-1)v]= 0$ for all $v$ in  $V$, and $\delta=1$ otherwise. If $t=n$ and  $\sigma$ is  in $O^-_{2n}(\F)$, then $\delta=1$.    If $t$ is odd, then $\delta = 1$.
\end{Rem}

For the last result in this section, $\ell = 3$.

\begin{proposition} \label{exc3}
Let $V$ be an irreducible symplectic $\F_3[G_\Q]$-module with squarefree conductor $N$.  Set $2n = \dim_\F V,$ $F= \Q(V)$ and $G = \Gal(F/\Q)$.   Then 
\begin{enumerate} [\rm i)]
\item $G \simeq \GS_{2n}(\F_3)$  or
\item $n$ is even, $G \simeq \SP_{n}(\F_3) \wr C_2$ and  $N=\gn \overline{\gn}$ in $\Z[\Mu_3]$. \end{enumerate}
\end{proposition}

\proof
An irreducible proper subgroup of $\SL_{2n}(\F_3)$ generated by transvections is isomorphic to $\SP_{2n}(\F_3)$, cf. \cite{KM}.  The pairing on $V$ implies that $F$ contains $\Mu_3$. The subgroup $P$ of $G$ generated by all transvections fixes $K =\Q(\Mu_3)$ and $F^P$ is unramified outside $3\infty$, so $F^P = K$ by Lemma \ref{Onlyell}.  If $V_{|P}$ is irreducible, then (i) holds.  If $V_{|P}$, is reducible, the arguments in the proofs of  Lemma \ref{Cliff} and Prop.\! \ref{exc2} give (ii), with $\gn$ a generator for the conductor ideal of the $\F[G_{\Q(\Mu_3)}]$-module $W$.   \qed

\section{Discriminants of stem fields}  \label{DiscSect}  
Let $F/k$ be a Galois extension of number fields with group $G$.  Let $\cD$ the decomposition group of a fixed prime $\pi_F$ of $F$ and $\cI_m$ the $m^{\rm th}$ ramification group (see \S \ref{FontSect}), with $\cI = \cI_0$ the inertia group.   For intermediate fields $L$, set $\pi_L=\pi_F\cap L$.  

\begin{theorem} \label{cosetdisc}
Let $G$ act transitively on $X$. If  $K$ is the fixed field of $G_x$ and   $\cI_m \backslash X$   is the set of $\cI_m$-orbits of $X$, then
$$
\ord_{\pi_k}(d_{K/k}) = \sum_{m \ge 0} \frac{1}{[\cI:\cI_m]} \: \left(|X|-| \cI_m \backslash X |\right).
$$
\end{theorem}

\proof
  If $H=G_x$ and $I$ is any subgroup of $G$, then $HgI \leftrightarrow Ig^{-1}x$ is a bijection between the set of double cosets $H\backslash G/I$ and the set of orbits $I \backslash X$.  Thus, 
\begin{equation} \label{coset1}
\sum_{HgI \,\in \, H\backslash G/I} [I:(I \cap H^g)] = [G:H],
\end{equation} 
where  $H^g = g^{-1}Hg$.  Suppose further that $J$ is a normal subgroup of $I$, so that $(I \cap H^g)J = I \cap H^gJ$ is a subgroup of $I$.  For each $g \in G$, we have 
\begin{equation} \label{coset2}
HgI = \bigsqcup Hgz_iJ,
\end{equation}
where $z_i$ runs over a set of representatives for the right cosets $I/(I \cap H^g)J$.   The isomorphism
$
(I \cap H^g)/(J \cap H^g) \simeq (I \cap H^g)J/J,
$ 
implies that
\begin{eqnarray} \label{coset3}
\nonumber  \sum_{HgI \, \in \, H\backslash G/I} \frac{|J \cap H^g|}{|I \cap H^g|} 
&=&  \sum_{HgJ \, \in \, H\backslash G/J} 
  \frac{1}{[I:(I \cap H^g)J]} \: \frac{|J \cap H^g|}{|I \cap H^g|}   \\
           &=& \sum_{HgJ \, \in \, H\backslash G/J}      
    \frac{1}{[I:J]} = \frac{|H\backslash G/J|}{[I\!:\!J]}.
\end{eqnarray}

\medskip

The ramification groups for $\pi_F$ inside $H$ are given by $\cI_m \cap H$ and the different ideal $\gD_{F/k}$ satisfies
$
\ord_{\pi_F}(\gD_{F/k}) = \sum_{m=0}^{\infty} (|\cI_m|-1).
$
By transitivity of differents, 
\begin{eqnarray} \label{diffKQ}
\ord_{\pi_K}(\gD_{K/k}) 
\nonumber  &=& \frac{1}{|\cI \cap H|}  \ord_{\pi_F}(\gD_{K/k}) \\ 
\nonumber  &=& \frac{1}{|\cI \cap H|}     
    \left(\ord_{\pi_F}(\gD_{F/k}) 
                    - \ord_{\pi_F}(\gD_{F/K})\right) \\                   
           &=& \sum_{m \ge 0} \frac{|\cI_m| - 
            |\cI_m \cap H|}{|\cI \cap H|}. 
\end{eqnarray}

Each prime of $K$ over $\pi_k$ has the form $g(\pi_F) \cap K$, corresponding to a unique double coset $Hg\cD$ in $H \backslash G/\cD$.  Since the decomposition and inertia groups of $g(\pi_F)$ inside $G$ are $g \cD g^{-1}$ and $g \cI g^{-1}$, the ramification and residue degrees of $g(\pi_F) \cap K$ over $\pi_k$ are given by 
\begin{equation} \label{f}
e(Hg\cD) = [\cI:(\cI \cap H^g)] \;\; {\rm and} \;\; f(Hg\cD) = [\cD:(\cD \cap H^g) \cI] .
\end{equation}
By conjugation,  (\ref{diffKQ}) implies that the exponent of $g(\pi_F) \cap K$ in $\gD_{K/k}$ is
\begin{equation} \label{x}
x(Hg\cD) = \sum_{m \ge 0} \frac{|\cI_m| - |\cI_m \cap H^g|}{|\cI \cap H^g|}.
\end{equation}
Moreover,
\begin{equation} \label{globalDiff}
 \ord_{\pi_k}(d_{K/k}) = \sum_{Hg\cD \, \in \,   
            H\backslash G/\cD} x(Hg\cD)f(Hg\cD).
\end{equation}

In view of (\ref{coset2}) and (\ref{f}), $Hg\cD$ is the disjoint union of $f(Hg\cD)$ distinct elements of $H\backslash G/\cI$.  By (\ref{globalDiff}) and (\ref{x}), we now have
\begin{eqnarray*}
\ord_{\pi_k}(d_{K/k})= \hspace{-10 pt} \sum_{Hg\cI \, \in \, H\backslash G/\cI} \hspace{-10 pt} x(Hg\cD)  =  \sum_{m \ge 0} \; \sum_{Hg\cI \, \in \, H\backslash G/\cI}  
           \frac{|\cI_m| - |\cI_m \cap H^g|}{|\cI \cap H^g|}.
\end{eqnarray*} 
But (\ref{coset1}) implies that 
\begin{eqnarray*}
\sum_{Hg\cI \, \in \, H\backslash G/\cI} \frac{|\cI_m|}{|\cI \cap H^g|} = \sum_{Hg\cI \, \in \, H\backslash G/\cI} \frac{[\cI\!:\!(\cI \cap H^g)]}{[\cI\!:\!\cI_m]}   
      = \frac{[ G:H]}{[\cI\!:\!\cI_m]} = \frac{[K:k]}{[\cI\!:\!\cI_m]},
\end{eqnarray*}
while (\ref{coset3}) with $J = \cI_m$ gives
$$
\sum_{Hg\cI \, \in \, H\backslash G/\cI} \frac{|\cI_m \cap H^g|}{|\cI \cap H^g|} = \frac{|H\backslash G/\cI_m|}{[\cI\!:\!\cI_m]}. 
$$
Substituting the last two identities in the previous double sum proves our claim.   \qed

\medskip

\begin{corollary} \label{tamestem}
If $\pi_k$ is tame in $F$, with ramification degree $|\cI(F/k)| = \ell$ prime, then 
$ 
\ord_{\pi_k}(d_{K/k})  = (1-\ell^{-1}) (|X| -  |X^\cI|).
$ 
\end{corollary}

\proof
Thm.\! \ref{cosetdisc} implies the claim, since $\cI_1$ is trivial and there are $|X^\cI|$ orbits of size 1, while the others have size $\ell$.  \qed

\medskip
We now apply these results to semistable Galois modules $V$ of conductor $N$. We write $F=\Q(V)$ and $G=\Gal(F/\Q).$

\begin{corollary} \label{pdiscSL} Let  $t=\ord_p(N)\ge 1$ and  $s=\dim_\F V$.  If $G$ acts transitively on $X =V-\{0\}$ and $K = F^{G_x}$, then $\ord_p(d_{K/\Q}) = (1-\ell^{-1})\left(q^s-q^{s-t}\right).$ 
\end{corollary}

\proof
Our claim follows from Cor. \ref{tamestem}, since $\dim V^\cI = s-t$ by (\ref{condexp}).  \qed
\bigskip

Now assume that  $\ell= 2$ and  $V$ is symplectic of dimension $2n.$  Let $K$ be the fixed field of $G_x,$ where $G$ acts transitively on $X,$ as below:
\begin{enumerate}[i)]
\item $G \simeq \cS_m={\rm Sym}(X)$ and $V$ is the representation in Remark \ref{SnRep}.
\item  $X = \Theta_{2n}^-$ or $X=\Theta_{2n}^- - \{\theta_0\}$, with $\theta_0$  fixed by $G$.  
\end{enumerate}

\begin{proposition} \label{pdisc} Let   $\cI_v=\langle\sigma\rangle\subseteq G$ be an inertia group at $v$ over $p\,|\,N.$ 
\begin{enumerate}[{\rm i)}]
\item If $G \simeq \cS_m$ and $\sigma$ is the product of $s$ disjoint transpositions, then $\ord_p(d_{K/\Q}) = s$ and $\ord_p(N) = \min(s,n).$  \vspace{2 pt}
\item If $G \simeq \SP_{2n}(\F)$ or $O^\pm_{2n}(\F)$, then $\ord_p(d_{K/\Q}) = \frac{1}{4}q^n(q^{n}-q^{n-t}- \delta),$  with $\delta$ as in {\rm Remark \ref{tdelta}}.

\end{enumerate}
\end{proposition}

\proof
i) Since $|X^{\cI_v}| = m-2s$, we have $\ord_p(d_{K/\Q}) = s$ by Cor. \ref{tamestem} and, by (\ref{condexp}), $\ord_p(N) = \dim_\F \, (\sigma-1)(V) = \min(s,n)$.

\hspace{15 pt} ii) We give a proof for  $t = 1$. Thus $\sigma$ is a transvection and we choose a symplectic basis for $V$ as in \S \ref{grpthy}, such that $\sigma = \tau_{[e_n]}$.  
For the even theta characteristic $\theta(x_1,\dots,x_{2n})= \sum_{j=1}^n x_j  x_{n+j}$, by (\ref{sx}) and (\ref{tQ}), we have 
$$
\sigma(\theta+a) = \theta + a + (1 + [a, e_n]) \, e_n.
$$
Thus, $\sigma$ fixes $\theta+a$ if and only if $[a, e_n] = 1$.  Let $V' =({\rm span}\{e_n,e_{2n}\})^\perp$ and $\theta'(y)= \sum_{j=1}^{n-1} y_jy_{n+j}$.  Assume $[a,e_n]=1$ and write $a = y+a_ne_n + e_{2n}$ with $y$ in $V'$.  In  $\F/\wp(\F)$, we have
$$
{\rm Arf}(\theta+a) = {\rm Arf}(\theta)+\theta(a) = a_n + \theta'(y).
$$   
Hence $\theta+a$ is in  $\Theta_{2n}^-$ precisely when one the following conditions holds: 
$$
{\rm (a) } \;  a_n \in \wp(\F) \; {\rm and} \; \theta'(y) \not\in \wp(\F) \quad {\rm or} \quad
{\rm (b) } \; a_n \not\in \wp(\F) \; {\rm and} \;  \theta'(y) \in \wp(\F).
$$
If $n=1$, only (b) applies, yielding $\frac{1}{2}q$ choices of $a$.  If $n \ge 2$, $y$ is in $\wp(\F)$ exactly when $\theta'+y$ is in $\Theta_{2n-2}^{+}$.  Hence there are $\frac{1}{2}q \, |\Theta_{2n-2}^-|$ choices of $a$ in case (a) and  $\frac{1}{2}q \, |\Theta_{2n-2}^+|$ choices in case (b).   But  $|\Theta_{2n-2}^+|+|\Theta_{2n-2}^-| = |V'| = q^{2n-2}$ and so $|(\Theta_{2n}^{-})^{\cI_v}| = \frac{1}{2}q^{2n-1}$.    \qed

\begin{Def}\label{augm}
A semistable Galois module $V$ is {\em ordinary at $2$} if it is symplectic and $\ga^2 \, V = 0$, where $\ga$ is the augmentation ideal in $\F[\cI_\lambda]$ for any $\lambda$ over 2 in $F$. 
\end{Def}

Let $V$ be the Galois module of a finite flat group scheme $\cV$ over $\Z_2$. Then $\cI_\lambda$ acts trivially on $V ^{m}$ and $V ^{et}.$ If the biconnected subquotient $\cV^b$ is trivial then $(\sigma-1)(\sigma'-1)(V ) = 0$ for all $\sigma, \sigma'$ in $\cI_\lambda$, whence $V$ is ordinary.   If $\cV^b \neq 0$, then $\cI_\lambda$ is not even  a 2-group.

\medskip

We next  treat the power of 2 in $d_{K/\Q}$ when $V$ is ordinary.  

\begin{lemma}  \label{orbitsize} 
We have $\ga \, V \subseteq Z \subseteq V^{\cI_\lambda}$ 
 for some  maximal isotropic subspace $Z$ of $V$.   If $H=G_\theta$ stabilizes an odd theta characteristic $\theta$, then 
$
| \cI_\lambda/(\cI_\lambda \cap H) | \le \frac{1}{2}q^n.
$
\end{lemma}

\proof
Set $\cI=\cI_{\lambda}.$ Since $\ga^2 \, V = 0$ and $\cI\subseteq \SP(V),$ we find $\ga \, V \subseteq V^\cI = (\ga \, V)^\perp$.  Thus, $\ga \, V$ is contained in a maximal isotropic space $Z$ and, by duality, $Z\subseteq V^\cI$.

If $\Gamma$ is the subgroup of $\SP_{2n}(\F)$ fixing both $Z$ and $V/Z$ pointwise, then we have $(g-1)(g'-1)(V) = 0$ for all $g,g'$ in $\Gamma$.  Hence $\psi(g) = (g-1)\theta$ defines a homomorphism $\Gamma \to V$.  In the notation of (\ref{sx}), $\Gamma$ is generated by the transvections $\tau_{[z]}$ with $z$ in $Z$.  Since we may identify $(\tau_{[z]}-1)\theta$ with $\sqrt{1+\theta(z)} \, z$, the homomorphism $\psi$ takes values in $Z$.  We next verify the exactness of the sequence
\begin{equation} \label{orbitsizeseq}
0 \to \Gamma \cap H \to \Gamma \stackrel{\psi}{\to} Z \stackrel{\theta}{\to} \F/\wp\F \to 0.
\end{equation}  
Since $Z$ is isotropic, $\theta$ is linear on $Z$ and $\theta$ is surjective because it is odd.  Clearly $\theta(\psi(\tau_{[z]}))$ is in $\wp \F$.  Conversely, if $\theta(z) = a^2 + a$ and $y=(1/\sqrt{a})z$, then $\psi(\tau_{[y]}) = z$.  This proves exactness around $Z$ and the rest is clear.  

Finally, $\cI \subseteq \Gamma$ and therefore 
$
|\cI/(\cI \cap H)| \le |\Gamma/(\Gamma \cap H)| = \frac{1}{2}q^n.
$
\qed

\medskip

\begin{proposition} \label{ordSpDisc}
If $V$ is ordinary at 2 and $G$ is transitive  on $\Theta_{2n}^-$ or $\Theta_{2n}^- - \{\theta_0\},$ then $\ord_2(d_{K/\Q}) \le (q^n-2)(q^n-1-\epsilon)$,  where $\epsilon = 0$ or $1,$ respectively.  
\end{proposition}

\proof
Since $\cI$ is a 2-group, $\cI_0 = \cI_1$.  The definition of the upper numbering (see \S \ref{FontSect}) and the bound on wild ramification (\ref{con}iii) imply that  $\cI_2 = 1$.  By Thm.\! \ref{cosetdisc}, $\ord_2(d_{K/\Q}) = 2(|X| - |\cI \backslash X|)$.

By Lemma \ref{orbitsize}, each $\cI$-orbit of $X$ has at most $\frac{1}{2}q^n$ elements and there are at least $2|\Theta_{2n}^-|/q^n = q^n-1$ orbits when $\epsilon = 0$, proving  the claim.

If $\epsilon = 1$, $\cI$ fixes $\theta_0$. The theta characteristic $\theta_0 + z$ is odd exactly if $\theta_0(z)$ is in $\wp\F$.  By (\ref{orbitsizeseq}), there are $\frac{1}{2}q^n$ such $z \in Z,$ giving at least $\frac{1}{2}q^n-1$ orbits of size 1 for $\cI$ acting on $X$.  The number of orbits not  accounted for is at least 
$$\frac{|X|- (\frac{1}{2}{q^n}-1)}{\frac{1}{2}q^n} = q^n-2$$ 
and so
$
|\cI \backslash X| \ge \frac{1}{2}q^n-1 +  (q^n-2) = \frac{3}{2} q^n-3.
$
Hence our claim. \qed

\medskip

\begin{proposition} \label{ordSmDisc}
If $V$ is ordinary and $G$ is a transitive subgroup of $\cS_m$, then $\ord_2(d_{K/\Q}) \le 2\lfloor m/2 \rfloor,$ unless $m=4$ or $8$, when $\ord_2(d_{K/\Q}) \le 3m/2$.
\end{proposition}

\proof
We find lower bounds for the number of $\cI$-orbits and apply Thm.\! \ref{cosetdisc}.  Since there is at least one orbit, our claims hold for $m \le 4$.  Assume $m \ge 5$ and refer to the explicit representation (\ref{SnRep}).  Let $y_{i,j} \in Y$ denote the vector with non-zero entries only in coordinates $i$ and $j$.  Write $\overline{y} \in V$ for the coset of $y \in Y$ when $m$ is even and $\overline{y} = y$ otherwise.

Suppose distinct letters $i,j$ lie in the same $\cI$-orbit.  If we can find a permutation $\sigma$ in $\cI$ such that $\sigma(i)=j$ and $\sigma(k) = k$, then 
$
\overline{y}_{i,j} = (\sigma-1)( \overline{y}_{i,k}) \in \ga \, V
$
is fixed by $\cI$.  It follows that $\tau(y_{i,j}) = y_{i,j}$ for all $\tau$ in $\cI$ and so $\{i,j\}$ is an $\cI$-orbit.

A larger orbit can exist only if $m = 2n+2$ is even and $\cI$ contains a product of $n+1$ disjoint transpositions, say 
$$
\sigma = (1,n+2)(2,n+3)\cdots(n+1,2n+2).
$$  
Treat subscripts modulo $2n+2$, fix $k$ and  consider $j \not\in \{k,k+n+1\}$.  Then 
$$
\overline{x}_j := \overline{y}_{j,j+n+1} -\overline{y}_{k,k+n+1} = (\sigma-1)(\overline{y}_{j,k}) \in \ga \, V
$$
is fixed by $\cI$.  If $m \ne 8$, $\overline{x}_j$ has a unique representative $x_j \in Y$ with exactly 4 non-zero entries and so $\tau(x_j) = x_j$ for all $\tau$ in $\cI$. Since
$$
\tau(k) \in \bigcap_{j \not\in \{k,k+n+1\}} \{j, j+n+1,k, k+n+1\} = \{k,k+n+1\},
$$
 $\{k,k+n+1\}$ is an $\cI$-orbit and every $\cI$-orbit has  2 elements.  If $m=8$, the $\cI$-orbits have size at most 4, giving the weaker bound.   \qed

\section{Stem field discriminant for $\Q(A[\gl])$ in a special case} \label{ParaSect}
In this section, $A_{/\Q}$ is $\go$-paramodular, with good reduction at $2 $ and $V = A[\gl]$ is irreducible for some   prime $\gl$ of $\go$ with residue field $\F_2$.  Any $\go$-polarization of $A$ has odd degree, since the kernel of the associated isogeny to $\widehat{A}$ intersects $A[\gl]$ trivially and thus $V = A[\gl]$ is symplectic for the Weil pairing.  Let $F = \Q(V)$ and $G  = \Gal(F/\Q)$.  The elements of $V$ correspond to differences $\theta_i-\theta_j$ of the 6 odd theta characteristics and we view $G$ as a subgroup of $\cS_6$, via its action on $\Theta^-$.  Irreducibility of $V$ implies that $G$ has an orbit $\Sigma \subseteq \Theta^-$ of size 5 or 6.  If $H=G_\theta$  stabilizes $\theta$ in  $\Sigma$, then $K = F^H$ is a stem field for $F$, with $[K:\Q] = |\Sigma|$.

The following local building blocks will appear.  Let $\tE = \Q_2(\Mu_3,\sqrt[3]{2})$ and let $X$ be the $G_{\Q_2}$-module, 2-dimensional over $\F_2$ with $\Q_2(X) = \tE$.  From  the exhaustive list \cite{JR1} of 2-adic fields of low degree, or by class field theory, there is a unique quartic extension $\tM/\Q_2$ whose Galois closure $\tL$ has non-trivial tame ramification, necessarily of degree 3.  Then $\tM/\Q_2$ is totally ramified,  $\ord_2(d_{\tM/\Q_2}) = 4$, $\tL \supset \tE$ and $\Gal(\tL/\Q_2) \simeq\cS_4$, with inertia subgroup $\cA_4$.  

\begin{proposition} \label{2-disc}
$\ord_2(d_{K/\Q}) \le  4$ {\em(}resp.  {\em 6)} if $[K:\Q] = 5$  {\em(}resp.  {\em 6)}.
\end{proposition}

\proof
If $V$ is ordinary at 2, the result follows from Prop.\! \ref{ordSpDisc} or \ref{ordSmDisc}.  Hence we suppose $F$ has non-trivial tame ramification over 2.  Among primes over 2 in $K$, choose $\lambda$ with maximal ramification degree  $e_\lambda(K)$ and consider all possibilities.
\begin{enumerate}[i)]
\item \, $e_\lambda(K) = 5$.   Then $(2)\cO_K =\lambda^5$ or $\lambda^5 \lambda'$, depending on whether $K$ is quintic or sextic, and $\ord_2(d_{K/\Q}) = 4$ by tame theory.
\item \,  $e_\lambda(K) = 3$.  If $K$ is quintic, the worst case occurs when $(2)\cO_K = \lambda^3 (\lambda')^2$ and then we have
$$
\ord_2(d_{K/\Q}) = \ord_2(d_{K_\lambda/\Q_2}) + \ord_2(d_{K_{\lambda'}/\Q_2}) = 2 + 2 = 4.
$$
Suppose $K$ is sextic.  If $(2) \cO_K = (\lambda \lambda')^3$, or $\lambda^3$ with residue degree $f_\lambda(K) = 2$, we have $\ord_2(d_{K/\Q}) = 4$.  In the remaining cases, at most one more prime $\lambda'$ over 2 ramifies in $K$, with $e_{\lambda'}(K) = 2$ and we conclude as for quintics.
\item \, $e_\lambda(K)=4$.  Then the completion $K_\lambda = \tM$.  If $[K:\Q] = 5$, the other prime over 2 in $K$ is unramified, but if $[K:\Q] = 6$, there may at worst be some $\lambda'$ with $e_{\lambda'}(K) = 2$.  Hence
$$
\ord_2(d_{K/\Q}) \le \left\{ \begin{array}{ll}
       4 & \mbox{ if  $\; [K:\Q] = 5$}, \\
       4+2=6 &\mbox{ if  $\; [K:\Q] = 6$}.  \end{array} \right.
$$
\item \, $e_K(\lambda) = 6$, so $[K\!:\!\Q] = 6$, $(2)\cO_K = \lambda^6$ and the inertia group $\cI$ of $\lambda$ acts transitively on $\Theta^-$. Since a non-zero fixed point for the action of $\cI$ on $V$ corresponds to a pair of theta characteristics preserved by $\cI$, contradicting transitivity, there are none. The tame ramification group $\cI/\cI_1$ is a cyclic subgroup of $\cS_6$ whose order is odd and a multiple of 3.  Hence $|\cI/\cI_1|=3$.  

\hspace{8 pt} Because $\cI_1$ is a non-trivial 2-group, normal in its decomposition group $\cD$, the fixed space $W = V^{\cI_1}$ is a non-zero $\cD$-module, properly contained in $V$.  Viewed as an $\cI/\cI_1$-module, $W$ is semisimple.  But $\cI/\cI_1$ has no non-zero fixed points on $W$, as they would be fixed points of $\cI$, so $\dim W = 2$ and $W \simeq X$.

\hspace{8 pt} Let $\cV = A[\gl]$, viewed as a finite flat group scheme over $\Z_2$.  The multiplicative component $\cV^m$ cannot have order 4, since $\cI$ is not a $2$-group, nor can it have order $2$, since $\cI$ has no non-trivial fixed points.  Hence $\cV^m = 0$ and $\cV$ is fully biconnected.  There is a subgroup scheme $\cW$ of $\cV$ with $\cD$-module $W$, and $\cV/\cW$ is biconnected, so its $\cD$-module also is isomorphic to $X$. 

\hspace{8 pt} Schoof \cite[Prop.\! 6.4]{Sch} showed that if $V$ is an extension of $X$ by $X$ as a $\cD$-module, then $\Q_2(V)$ is contained in the maximal elementary 2-extension $\tL_1$ of $\tE$ with ray class conductor exponent 2.  One checks that $\tL_1$ is an unramified central extension of degree 2 over $\tL$ and the root discriminant of $\tL_1/\Q_2$ is $7/6$.  Since $\ord_2(d_{K/\Q})$ is even, we have $\ord_2(d_{K/\Q}) \le 6$, as claimed.  \qed
\end{enumerate}

\section{Preserving the Fontaine bound} \label{FontSect}  
Let $K'/K$ be a Galois extension of $\ell$-adic fields with Galois group $G$.  Denote the ring of integers of $K'$ by $\cO'$ and a prime element by $\lambda'$.  Set
$$
G_n = \{ \sigma \in G \, | \,
         \ord_{\lambda'}(\sigma(x)-x)  \ge n+1 \;  {\rm for \; all} \; x \in \cO' \},
$$
so $G_0$ is the inertia group and $t_{K'/K} = [G_0\!:\!G_1]$ is the degree of tame ramification.  If  $\lfloor x \rfloor = m$, the Herbrand function is given by
\begin{equation} \label{Herb} 
\varphi_{K'/K}(x) =\frac{1}{|G_0|}(\,|G_1|+\dots+|G_m|+(x-m)|G_{m+1}|\,)
\end{equation} 
and is continuous and increasing.  In the upper numbering of Serre \cite[IV]{Ser1}, $G^{m} = G_n$, with $m = \varphi_{K'/K}(n)$.  In the numbering of \cite{Fo} or \cite{JR1}, this group is denoted $G^{(m+1)}$.  Let $\psi_{K'/K}$ be the inverse of $\varphi_{K'/K}$.

\begin{Not} \label{Defcm}
Let $c = c_{K'/K} $ be the maximal integer such that $G_c \ne 1$.  We omit the lower field if $K=\Q_{\ell}$.  Let $m_{K'} = \psi_{K'/\Q_\ell}(\frac{1}{\ell-1})$. 
\end{Not}

Wild ramification in $K'/K$ is equivalent to $c_{K'/K} \ge 1$.  If $G_1$ is not abelian, then $c_{K'/K} \ge 2$, since successive quotients in the ramification filtration are elementary abelian $\ell$-groups.  By (\ref{Herb}), $m_{K'}$ is an integer when $(\ell-1)$ divides $t_{K'/\Q_\ell}$.

\begin{lemma}   \label{buildFont}
Let $E \supset F$, both Galois over $K$, $G = \Gal(E/K)$ and $H = \Gal(E/F)$.  Then 
$
1 \to H^{\psi_{F/K}(x)} \to G^x \stackrel{res}{\longrightarrow} \Gal(F/K)^x \to 1
$
is exact.  In addition, 
\begin{equation} \label{ct}
m_E \ge t_{E/F} m_F \quad {\rm and} \quad  c_{E/K} \ge t_{E/F} c_{F/K}.
\end{equation}
\end{lemma}

\proof
By compatibility with quotients, {\em res} is surjective and its kernel is
$$
G^x \cap H = G_{\psi_{E/K}(x)} \cap H = H_{\psi_{E/K}(x)} = H^{\varphi_{E/F}\psi_{E/K}(x)} = H^{\psi_{F/K}(x)},
$$
since $\psi_{E/K} = \psi_{E/F} \, \psi_{F/K}$.  Thus the sequence is exact.

Def.\! \ref{Herb} implies that $ t_{E/F}\varphi_{E/F}(z) \le z$, so $\psi_{E/F}(z) \ge  t_{E/F}z$.   If $x=\varphi_{F/K}(c_{F/K})$, then $G_{\psi_{E/K}(x)} = G^x\neq 1$ by surjectivity of {\em res}.  Hence
$$
c_{E/K} \ge \psi_{E/K}(x) = \psi_{E/F}  \psi_{F/K}(x) = \psi_{E/F}(c_{F/K}) \ge  t_{E/F} \, c_{F/K}
$$
and similarly for $m_E \ge  t_{E/F}\, m_F$.   \qed

\begin{Def}  \label{FontDef}
Let $F$ be the Galois closure of $K/\Q_\ell$.  We say $K$ is {\em Fontaine} if $\Gal(F/\Q_\ell)^u = 1$ for all $u > \frac{1}{\ell-1}$, or equivalently, $c_F \le m_F$.
\end{Def}

\begin{lemma} \label{numbers} 
Let $E \supset F$, both Galois over $\Q_\ell$, $G = \Gal(E/\Q_\ell)$ and $H = \Gal(E/F)$.
\begin{enumerate}[{\rm i)}]
\item If $t_{F/\Q_\ell} = \ell-1$, then $m_E \ge t_{E/F}$, with equality when $G_0$ is abelian.
\item Let $F$ be Fontaine, with non-trivial wild ramification.  Then $1 \le c_F \le m_F$.  Assume further that  $t_{F/\Q_\ell} = \ell-1$.  Then $c_F = m_F=1$ and, if $E$ is Fontaine, then $c_E = m_E$.
\end{enumerate}
\end{lemma}

\proof 
i) Since $\varphi_{F/\Q_\ell}(1)=\frac{1}{\ell-1}$, we have $m_F = 1$, so $m_E \ge t_{E/F}$ by (\ref{ct}).  If $G_0$ is abelian and $t_{E/\Q_\ell}$ does not divide $j$, then $G_j = G_{j+1}$ by \cite[IV,\S 2]{Ser1}.  Thus the definition gives $\varphi_{E/\Q_\ell}(t_{E/F}) = \frac{1}{\ell-1}$, whence $m_E =  t_{E/F}$.

\vspace{2 pt}

ii) By Def.\! \ref{FontDef}, $\varphi_{F/\Q_\ell}(c_F) \le \frac{1}{\ell-1} = \varphi_{F/\Q_\ell}(m_F)$.  Hence $c_F \le m_F$.  If $t_{F/\Q_\ell} = \ell-1$, then $m_F = 1$, so $c_F = 1$.  Surjectivity of $res$ in Lemma \ref{buildFont} implies that $G^{\frac{1}{\ell-1}} \ne 1$.  If $E$ is Fontaine, it follows that $c_E = \psi_{E/\Q_\ell}(\frac{1}{\ell-1}) = m_E$.  \qed

\begin{Ex} \label{S4Q2} 
By class field theory or the table of quartics \cite{JR1}, there is a unique Fontaine $\cS_4$-extension $F/\Q_2$.  The ramification subgroups of $\overline{G} = \Gal(F/\Q_2)$ are $\overline{G}_0 \simeq \cA_4$, $\overline{G}_1 \simeq C_2^2$ and $\overline{G}_2 = 1$, so $c_F = 1$, $\varphi_{F/\Q_2}(x) = (4+(x-1))/12$ if $x \ge 1$ and $m_F = 9$.  Moreover, $E = F(i)$ remains Fontaine, with $G = \Gal(E/\Q_2) \simeq \cS_4 \times C_2$.  Lemma \ref{numbers}ii may be used to show that $\vert G_0 \vert = 24$, $\vert G_1 \vert = 8$, $\vert G_2 \vert = \dots = \vert G_9 \vert = 2$, $\vert G_{10} \vert = 1$ and $c_E = m_E = 9$.  Alternatively, $E$ has two stem fields of degree 6 and this determines $E$ uniquely in \cite{JR1}.
\end{Ex}

\begin{lemma} \label{raybound2} 
Let $M/F$ be abelian, with $F/\Q_\ell$ Galois.  Then $M$ is Fontaine if and only if $F$ is Fontaine and the ray class conductor exponent $\gf(M/F)  \le \lfloor m_F \rfloor +1$.  
\end{lemma} 

\proof 
If $E$ is the Galois closure of $M/\Q_\ell$,  then $E/F$ is abelian and we have $\gf(E/F) = \gf(M/F) = \varphi_{E/F}(c_{E/F})+1$, cf. \cite[XV, \S2]{Ser1}.   The exact sequence of Lemma \ref{buildFont} with $K = \Q_\ell$ implies our claim. \qed

\begin{Rem} \label{betterFont}
Let $E$ be a {\em number field} with root discriminant $\varrho_E$.  Write $\tE$ for the completion of $E$ at a prime $\lambda\,|\,\ell$ and $e_{\tE}$ for the absolute ramification degree.  Suppose $E$ contains $F$, both Galois over $\Q$, with $\tE$ Fontaine.  Then
\begin{center}
$ 
\ord_\ell(\varrho_{E}) \le 1 + \frac{1}{\ell-1} - \frac{t_{\tE/\tF} \, c_{\tF} +1}{e_{\tE}}. 
$ 
\end{center}
Indeed, if $\gD_{\tE/\Q_\ell}$ is the different, then \cite[IV, Prop.\! 4]{Ser1} and \cite[Prop.\! 1.3]{Fo} give
\begin{center}
$
\ord_\ell(\varrho_{E}) = \frac{1}{e_{\tE}} \, \ord_\lambda(\gD_{\tE/\Q_\ell}) = 1 + \varphi_{\tE/\Q_\ell}(c_{\tE}) - \frac{c_{\tE}+1}{e_{\tE}}.
$
\end{center}
Conclude by Def.\! \ref{FontDef} and (\ref{ct}). \qed
\end{Rem}

Because the upper numbering is compatible with quotients, the composition of Fontaine fields is Fontaine and there is a maximal field $L$, such that $\Gal(F/\Q_\ell)^u = 1$ for all Galois subfields $F$ finite over $\Q_\ell$ and all $u > \frac{1}{\ell-1}$.  Since $F = \Q_\ell(\Mu_\ell, (1-\ell)^ \frac{1}{\ell})$ is contained in $L$, Lemma \ref{numbers}ii implies a gap in the upper numbering:
$$
\Gal(L/\Q_\ell)^ \frac{1}{\ell-1} \ne \Gal(L/\Q_\ell)^{ \frac{1}{\ell-1}+\epsilon} \; {\rm for } \; \epsilon > 0.
$$

Hajir and Maier \cite{HM} study number field extensions $K'/K$ of {\em bounded depth}, i.e.\! with vanishing ramification groups $\cD_\gp(K'/K)^x$ for all $x \ge \nu_\gp.$  When there is deep wild ramification, the concept of Galois slope content introduced by Jones and Roberts \cite{JR1} and used in \cite[\S 1.4]{JJ}, leads to variants of (\ref{ct}) and Remark \ref{betterFont}, not required for our applications, thanks to (\ref{con}iii).

\section{Using Odlyzko} \label{Odl}
We study some maximal $(\ell,N)$-controlled extensions $L/\Q$ by means of Odlyzko's bounds \cite{Od,Od2,dyd}.  If the $\F[G_\Q]$-module $V$ is semistable and bad only at $S$, then $\Q(V)$ is $(\ell,N_S)$-controlled.  The converse holds for $\ell=2$  but not for $\ell$ odd; e.g.\!  if $\dim V=2$, then ${\rm Sym}^2V$ rarely is semistable.

By tameness at $p \, \vert \, N$ and the bound (\ref{con}iii), the root discriminant of $L/\Q$ satisfies $\varrho_L < \ell^{1+\frac{1}{\ell-1} } \, N^{1-\frac{1}{\ell}}$.  More precisely:
\begin{equation} \label{GFdisc}
\ord_p(\varrho_L) \le 1-\ell^{-1} \;{\rm for \; all } \; p \,\vert\, N \quad {\rm and} \quad \ord_\ell(\varrho_L)<1+(\ell-1)^{-1}.
\end{equation}
 
\begin{proposition} \label{Onlyell}  
For $\ell\le 13$, the maximal $(\ell,1)$-controlled extension $L$ is $\Q(\mu_{2 \ell})$.  Under GRH, the same is true for $\ell=17$ and $19$.
\end{proposition}

\proof  
For $\ell$ odd, $\Q(\Mu_{\ell}) \subseteq L$ and $n = [L\!:\!\Q]$ is a multiple of $\ell-1$.  From (\ref{GFdisc}) and \cite{Od}, we find $M$ in Table \ref{ell1} below  such that $n \le (\ell-1)M$.  If $\ell = 13, 17, 19$, we see that $M < \ell$, so $L/\Q$ is tame at $\ell$ and $\varrho_L \le  \ell^{1-\alpha}$, with $\alpha=((\ell-1)M)^{-1}$.  One gets a new bound $n \le (\ell-1)M'$ with $M' \le 5$.  If $\ell \le 11$, we have $M \le 5$.  In both cases,  $L$ is abelian over $\Q(\Mu_{2\ell})$ and so $L = \Q(\Mu_{2\ell})$ by class field theory \cite[Lem.\! 2.2]{BK1}.  Use $\Q(i) \subseteq L$ for $\ell = 2$.  \qed

\vspace{-10 pt}

{\small
\begin{table}[here]
\begin{center}
\begin{tabular}{|c||c|c|c|c|c|c|c|c|}
\hline
 $\ell$ &2 & 3 & 5 &  7 & 11 & 13 & 17 &19 \\
\hline
$\varrho_L\le$ & 4 & 5.197& 7.477& 9.682&13.981& 16.099& 20.294& 22.377  \\
\hline
 $M$ & 2 & 3  & 3  &  3 & 5 & 7 &$ 8$ & $ 10$ \\
\hline
\end{tabular} 
\smallskip
\begin{caption}{Odlyzko bounds for $(\ell,1)$-controlled fields} \label{ell1}\end{caption}
\end{center}
\end{table}}

\vspace{-10 pt}

Now suppose $L$ is maximal $(2,N)$-controlled, so $\varrho_L < 4 N^\frac{1}{2}$ by (\ref{GFdisc}).  If $n = [L\!:\!\Q]$ is finite, \cite[Tables 3,4]{Od} provides $B, E$, depending on a parameter $b$, such that
$
\varrho_L > B e^{-\frac{E}{n}}.   
$
In Table \ref{2Nbd} below, we find a best  bound for $n < E/\log({B}/{4 N^\frac{1}{2}})$ by varying $B > 4 N^\frac{1}{2}$, unconditionally for $N \le 21$ and under GRH for larger $N$. 

If $V$ is an irreducible semistable $\F_2[G_\Q]$-modules good outside $S$ and $N_S \vert N$, then $\Gal(\Q(V)/\Q)$ factors through $\overline{G}=G/H$, where $H$ is the maximal normal 2-subgroup of $G=\Gal(L/\Q)$.  For odd $N \le 79$ and $N =97$, we find a subfield $F$ of $L$ containing $L^H$ by composing a solvable extension of $\Q$ with a subfield of $\Q(J_0(N)[2])$.  Then we use the improvements in \S \ref{FontSect} on the bound (\ref{GFdisc}) for $\varrho_L$, together with the Odlyzko tables and Magma \cite{BCP} to control $[L\!:\!F]$. 

{\small
%\begin{center}
\begin{table}[here]
\begin{tabular}{|c||c|c|c|c|c|c|c|c|c|c|c|c|c|}
\hline
 $N$ & 3 & 5 & 7 &11& 13 &15&17&19 &21& 23& 29\\
\hline {$n\le$}  & 10  &
  16 & 22&42 & 56  & 74 & 100&138 &192 &  98&155\\
\hline
%\hline
 $N$ & 31& 33 &35& 37& 39& 41& 43&47& 51&53&55\\
\hline {$n\le$}&181&210 &244&284&330&385&449&615&  852&1007&1196\\
\hline
%\hline
 $N$&57&59&61  &65&67&69&71& 73&77&79&97\\ 
\hline {$n\le$}&1427&1710&2061 &3046&3743&4638&5800&  7332&12042&15766 & 470652\\
\hline
\end{tabular}
\smallskip 
\begin{caption}{Bounds on $n = [L\!:\!\Q]$ for $(2,N)$-controlled fields $L$} \label{2Nbd} \end{caption}
\end{table}
%\end{center} 
}

\begin{theorem}[GRH] \label{Odl2} Let $V$ be semistable and irreducible over $\F_2$.  If $V$ is bad exactly over $S$ and $N = N_S \le 79$ or $N=97$, the following hold.  
\begin{enumerate}[{\rm i)}]
\item  No such $V$ exists for $N$ in $\{3, 5, 7, 13, 15, 17, 21, 33, 39,41, 55, 57, 65, 77\}.$   
\item  $V$ is unique and $\dim V = 2$ for $N$ in $\{11,19, 23, 29, 31, 35, 37, 43, 51, 53, 61\}.$
\item  $V$ is unique for $N$ in $\{23,31,47, 71\}$.
\item $V$ is an irreducible $\F_2[\overline{G}]$-module with $\overline{G} = D_9$, $D_3 \times \cA_5$, $\cA_5$, $D_3 \times D_5$, $\SL_2(\F_8)$ when $N = 59, \, 67, \, 73, \, 79, \, 97$ respectively.
\end{enumerate}
\end{theorem} 

\begin{Rem}  
Aside from $\F_2$, there are exactly two irreducible $\F_2[\cA_5]$-modules, both 4-dimensional, occurring as a submodule $V_1$ and quotient module $V_2$ of the permutation module. The non-trivial $\F_2[\SL_2(\F_8)]$-modules have dimensions 6, 8 and 12. Further,  the irreducible modules for $G_1 \times G_2$  are the tensor products of irreducibles for $G_1$ and $G_2$.
\end{Rem}

\noindent{\em Sketch of Proof.} In (i), $G$ is a 2-group, except for 33, 55, 57, 77, when  $\overline{G} \simeq D_3$ has a representation whose conductor, 11 or 19, divides $N$ properly.  In (ii), $V \simeq C_N[2]$ for an elliptic curve $C_N$ of conductor $N$, except that $V \simeq J_0(29)[\sqrt{2}]$ for $N=29.$ In (iii), $V$ is the $\F_2[D_h]$-module of dimension $h-1$ induced by the Hilbert class field over $\Q(\sqrt{-N})$ of class number $h=3,$ 3, 5, 7 corresponding to $N = 23$, 31, 47, 71 respectively.  

${\mathbf N=59\!:} \;$ The two irreducibles are the constituents of $J_0(59)[2]$, using an equation for $X_0(59)$, namely $y^2 = f(x)g(x)$ with $f = x ^3 - x^2-x + 2$ and 
$$g = x^9 -7 x^8 +16 x^7-21 x^6+12x^5 - x^4-9 x^3+6 x^2-4 x-4.$$
The Galois group of $g$ is $D_9$ and a root of $f$ gives a cubic subfield. 

${\mathbf N=67\!:} \;$ Let $V_1=C_{67}[2]$  and  $V_2=J_0^+(67)[2].$ Then $\Gal(\Q(V_2)/\Q)=\SL_2(\F_4)$  and $[L\!:\!\Q(V_1,V_2 , i)]\le 2.$ 

We provide more details for $N=73$, 77, 79 and 97.  Let $E$ be the maximal abelian  extension of $\Q$ in  $L$.  Since $G$ is generated by involutions,  $E/\Q$ is the elementary 2-extension generated by $i$ and $\sqrt{p}$ as $p$ ranges over $S$.   

\begin{lemma} \label{notame}
Let $M \supset F$ be subfields of $L$ containing $E$ and Galois over $\Q$.  Set $T = \Gal(M/F)$ and assume $\lambda \, \vert \, 2$ is totally ramified of odd degree $t = \vert T \vert > 1$ in $M/F$.  Then $t=3$ and the residue degree $f_\lambda(E/\Q)=2$.
\end{lemma}

\proof
Since the image of $\alpha\!:  \Gal(M/\Q) \to \Aut(T)$ by conjugation is abelian, $E$ contains $M_0 = L^{\ker \alpha}$ and so $f = f_\lambda(M_0/\Q) \le f_\lambda(E/\Q) \le 2.$ Any Frobenius in $\cD_\lambda(M/M_0)$ acts trivially on $T.$ Thus $2^f \equiv 1 \! \pmod{t}$ and the claim ensues.  \qed 

\medskip

\begin{Rem} \label{boot}
Let $M \supseteq F$ be subfields of $L$ containing $\Q(i,\sqrt{N})$ and Galois over $\Q$.  Denote the residue, ramification and tame degree of $\lambda$ in $F/\Q$ by  $f_0$, $e_0$ and $t_0$ respectively.  Given an {\em a priori} bound $[M\!:\!F] \le b$, consider possible factorizations $[M\!:\!F] = 2^s  t_1 u_1$, where $2^s$ is the degree of wild ramification, $t_1$ the degree of tame ramification and $u_1=f_1g_1$ the unramifed (inert and split) degree of $\lambda$ in $M/F$.  The resulting tame ramification in $M/\Q$ requires that the completion $M_\lambda$ contain $\Mu_{t_0t_1}$ and so $2^{f_0f_1} \equiv 1 \pmod{t_0t_1}$.  

For each $s$ with $0 \le s \le \log_2 b$, let $t_1 \ge 1$ run through odd integers at most $b/2^s$.  Set $\beta = (c_F \, t_1+1)/(2^s \, t_1 \, e_0)$, as in Remark \ref{betterFont} and let $n_\beta$ be the Odlyzko bound on $[M\!:\!\Q]$ when $\varrho_M \le 2^{2-\beta}\sqrt{N}$.  Then 
$
1 \le g_1 \le {n_\beta}/({2^s t_1 f_1 \, [F\!:\!\Q]}).
$
Values of $s,t_1,f_1$ not satisfying the congruence and inequality above are ruled out.
\end{Rem}
  
Let $E_1$ be the maximal subfield of $L$ abelian over $E$.  By Lemmas \ref{numbers}ii and \ref{raybound2}, the ray class conductor of $E_1/E$ divides $(1+i)^2 \cO_{E}$.  Then class field theory or Magma gives Table \ref{7397} below.

{\small \begin{table}[here]
\begin{center}
\begin{tabular}{|c|c|c|c|c|}
\hline
$N$ & $\Gal(E_1/E)$ & $e_\lambda(E_1/\Q)$ & $f_\lambda(E_1/\Q)$ & $g_\lambda(E_1/\Q)$ \\ \hline
73 &  $C_4$      &  4 & 2 & 2 \\ \hline
77 &  $C_6$      &  6 & 2 & 4 \\ \hline
79 &  $C_{15}$  &  2 & 5 & 6 \\ \hline
97 & $C_4$        & 4 & 2 & 2  \\ \hline
\end{tabular}
\smallskip
\begin{caption}{Decomposition type of $\lambda \vert 2$ in $E_1$} \label{7397} \end{caption}
\end{center}
\end{table}}

\vspace{-15 pt}

${\mathbf N = 73\!:} \;$  The Jacobian $J_0^+(73)$ has RM by $\Q(\sqrt{5})$ and the Galois group of its 2-division field $K$ is $\SL_2(\F_4) \simeq \cA_5$.  For the 5 primes over 2 in $K$, $f_\lambda(K/\Q) = 3$ and Frobenius acts irreducibly on $\cI_\lambda(K/\Q) \simeq C_2^2$.  Since Frobenius is reducible on $\cI_\lambda(E_1/\Q) \simeq C_2^2$, we have $\cI_\lambda(F/\Q) \simeq C_2^4$ for the compositum $F = E_1K$, thus $[F\!:\!\Q] = 960$.  By Table \ref{2Nbd}, $[L\!:\!\Q]=960r \le 7332$, so $r \le  7$.  Lemma \ref{notame} implies the tame degree $t_\lambda(L/F) = 1$, so $e_\lambda(L/F)$ divides 4. Finally, $[L\!:\!F]$ divides 4 by (\ref{boot}).

${\mathbf N = 77\!:} \;$  In the $\cS_3$-field $K_0 = \Q(J_0(11)[2]) = \Q(\sqrt{-11},\theta)$, with $\theta^3-2\theta^2+2 = 0$, the decomposition type over 2 is $e_\lambda = 3$, $f_\lambda = 2$, $g_\lambda = 1$.  If $K = E(\theta) = K_0(i, \sqrt{-7})$, then $\Gal(K/\Q) \simeq C_2 \times C_2 \times \cS_3$ and $\cI_\lambda(K/\Q) \simeq C_6$, so $m_K = 3$ by (\ref{numbers}i).   If $F$ is the maximal subfield of $L$ abelian over $K$, the ray class conductor of $F/K$ divides $(1+i)^4 \cO_K = 4\cO_K$ by Lemma \ref{raybound2}.  Then $\Gal(F/K) \simeq C_2 \times C_2 \times C_4$ and the decomposition type of 2 is
$
e_\lambda(F/\Q) = 48, \,  f_\lambda(F/\Q) = 2 \; {\rm and } \;  g_\lambda(F/\Q) = 4.
$  

A group of order $3 \cdot 2^a$ admits a unique quotient isomorphic to $C_3$ or $\cS_3$. If $[L\!:\!K_0]=3 \cdot 2^a$, there is a  $C_3$ or $\cS_3$ extension of $K_0$. The latter provides a central quadratic $M_0/K_0$, with $M_0/\Q$ Galois and $\Gal(M_0/\Q) \simeq D_6$. In both cases, we find that $\Gal(M_0K/K) \simeq C_3$, contradicting $[F\!:\!K] = 16$.

We claim that $\Gal(L/F)$ is a 2-group.  If not, since $[L\!:\!F] \le 31$ from Table \ref{2Nbd} and  $[L\!:\!F] \ne 3 \cdot 2^a$, the wild ramification degree $\vert \cI_\lambda(L/F)_1 \vert$ divides 4.  Example \ref{S4Q2} and (\ref{ct}) imply that $\cI_\lambda(F/\Q)_9 \ne 1$.  Use Remark \ref{boot} with $c_F \ge 9$ to show that the only remaining case is $[L\!:\!F] = 10$, with tame degree $t_\lambda(L/F) = 5$ and wild degree 2. It is precluded by Lemma \ref{notame}.

Thus the kernel of the surjection $G \stackrel{\eta}{\to} \Gal(\Q(J_0(11)[2])/\Q) \simeq \cS_3$ is a 2-group and irreducible representations $V$ of $G$ factor through $\Image \eta$, of conductor 11, so there is no $V$ of conductor $77$.

${\mathbf N = 79\!:} \;$  The strict class fields $H^\pm$ of $\Q(\sqrt{\pm 79})$ have respective orders 3 and 5 and so $E_1 = H^+H^-.$  Let $K^\pm$ be the maximal subfields of $L$ abelian respectively over $H^\pm(i)$.  Since $e_\lambda(H^\pm(i)/\Q) = 2$, the ray class conductors of $K^{\pm}/H^\pm(i)$ divide $(1+i)^2 \cO_{H^\pm(i)}$ by Lemma \ref{raybound2}.  Magma provides the following information.
\begin{eqnarray*}
&\Gal(K^+/H^+(i))& \simeq C_2^2 \times C_3, \hspace{3 pt} \; {\rm with } \; e_\lambda = \, \, 2, \, f_\lambda= 2, \, g_\lambda = 3.\\
&\Gal(K^-/H^ -(i))& \simeq C_2^4 \times C_5,  \hspace{3 pt} \; {\rm with } \; e_\lambda = 16, \, f_\lambda= 5, \, g_\lambda = 1. 
\end{eqnarray*}

If $\Gal(E_1/E) = \langle \tau \rangle$, then $\tau^5$ and $\tau^3$ act trivially on $K^+$ and $K^-$ respectively.  Hence $\tau$ is trivial on $K^+ \cap K^-$ and $(K^+ \cap K^-)/E$ is abelian.  Since $K^+ \cap K^-$ contains $E_1$, equality holds by maximality of $E_1$.  For $F = K^+K^-$, we therefore have $[F\!:\!E_1] = 2^6$,  $[F\!:\!\Q] = 3840$ and $[L\!:\!\Q] = 3840r$, with $r \le 3$.  Because Frobenius acts irreducibly on $\cI_\lambda(K^-/E_1) \simeq C_2^4$ but trivially on $\cI_\lambda(K^+/E_1) \simeq C_2$, we see that $\cI_\lambda(F/E_1) \simeq C_2^5$ and $e_\lambda(F/\Q) = 64$.  By Lemma \ref{notame}, $t_\lambda(L/F) = 1$, so $[L:F] \le 2$ by (\ref{boot}).  Thus the kernel of $G \twoheadrightarrow \Gal(H^+/\Q) \times \Gal(H^-/\Q) \simeq D_3 \times D_5$ is a 2-group.

${\mathbf N = 97\!:} \;$  There is a subfield $K$ of $\Q(J_0(97)[2])$ with $\Gal(K/\Q) \simeq \SL_2(\F_8)$ and decomposition type 
$
e_\lambda(K/\Q) = 8, \, f_\lambda(K/\Q) = 7, \, g_\lambda(K/\Q) = 9.
$
Any Frobenius in $\cD_\lambda(K/\Q)$ acts irreducibly on $\cI_\lambda(K/\Q) \simeq C_2^3$ but reducibly on $\cI_\lambda(E_1/\Q) \simeq C_2^2$, so $\cI_\lambda(F/\Q) \simeq C_2^5$ for the compositum $F = E_1K$.   Since $[F\!:\!\Q] = 504 \cdot 16 = 8064$, Table \ref{2Nbd} implies that $[L:F] \le 58$.  Thus the dimensions of irreducible representations of $\SL_2(\F_8)$ over $\F_p$ for small $p$ force the action of $\Gal(F/E_1)$ on the maximal abelian quotient of $\Gal(L/F)$ to be  trivial.  But no central extension of $\SL_2(\F_8) $ is perfect  \cite{Asch,HBI}. Hence $L$ is the compositum of $F$ with a solvable extension of $E_1$.  The ray class extension of $E_1$ whose conductor divides $\prod \lambda^2$, as $\lambda$ runs over the primes above 2 in $\cO_{E_1}$, turns out to be trivial, whence $L=F$ by Lemmas \ref{numbers}ii and \ref{raybound2}.   \qed

\medskip

${\mathbf N = 127\!:} \;$  We begin the solvable tower with $E_0 = E = \Q(i,\sqrt{127})$ and find successive maximal abelian extensions $E_{j+1}/E_j$ in $L/\Q$.  For ray class conductor $(1+i)^2 \cO_E$, we have $[E_1\!:\!E] = 5$.  Thus $E_1$ is the compositum of $\Q(i)$ and the Hilbert class field over $\Q(\sqrt{-127})$.  Now $e_\lambda(E_1/\Q) = 2$, so the ray class conductor of  $E_2/E_1$ divides $(1+i)^2 \cO_{E_1}$ and we have $\Gal(E_2/E_1) = C_2^4$.   Moreover, any Frobenius in $\cD_\lambda(E_2/\Q)$ has irreducible action of order 5 on this ray class group.  The decomposition type over 2 is 
$
e_\lambda(E_2/\Q) = 32, \, f_\lambda(E_2/\Q) = 5, \, g_\lambda(E_2/\Q) = 2.
$
The ray class conductor of $E_3/E_2$ divides $\prod \lambda^2$, as $\lambda$ runs over the primes of $\cO_{E_2}$ above 2, but we do not know whether $E_3 = E_2$.  

There is a subfield $K$ of $\Q(J_0(127)[2])$ with $\Gal(K/\Q) \simeq \SL_2(\F_8)$ and decomposition type 
$
e_\lambda(K/\Q) = 8, \, f_\lambda(K/\Q) = 7, \, g_\lambda(K/\Q) = 9.
$
Any Frobenius in $\cD_\lambda(K/\Q)$ has irreducible action of order 7 on $\cI_\lambda(K/\Q) \simeq C_2^3$.  For the compositum $F = E_2K$, of degree $320 \cdot 504 = 161280$, we therefore have $\cI_\lambda(F/\Q) \simeq C_2^8$ and so $c_F = m_F = 1$ by Lemma \ref{numbers}.  By Remark \ref{betterFont}, the root discriminant is $\varrho_F = 2^{2-\frac{1}{128}} \sqrt{127} \approx 44.834.$  Since this just exceeds the asymptotic bound $8\pi e^{\gamma}\approx 44.763$, where $\gamma$ is Euler's constant \cite{Od2}, we do not know whether $[L:\Q]$ is finite and thus it would be entertaining to find $L$.

\section{Some speculations} \label{spec}
Assume $\ell \nmid N$ is prime. Let AB$(N)$  be the set of simple semistable abelian varieties $A$ with   $N_A\vert N$ and let Irr$(\ell,N)$ be the set of irreducible semistable $\F_\ell[G_\Q]$-modules $V$ with  $N_V\vert N$, both taken up to isomorphism.

\vspace{5 pt}

{\bf Q1}. Is AB$(N)$ finite? 

{\bf Q2}. Is Irr$(\ell,N)$ finite? 

{\bf Q3}. Is $\log N_A>>g \log g$ as  $g=\dim A\to\infty$?

\vspace{5 pt}

\noindent Faltings' theorem  answers {\bf Q1} in the affirmative if  $\dim A$ is bounded  in terms of $N_A$.   Mestre \cite{Mes} shows that   standard {\em conjectures} on Hasse-Weil L-functions   imply that   $\dim A\le .42\, \log N_A,$ without assuming semistability or simplicity.  Such strong modularity or even an answer to  {\bf Q2} are  equally unavailable. However,  the negation of  {\bf Q1} violates   one of the following two {\em plausible} assertions.

\vspace{5 pt}

\setlength{\hangindent}{58 pt}
{\bf P$_1(\ell,N)$}. The dimension of the composition factors of $A[\ell]$ is  bounded for  $A$  in AB$(N)$.  

{\bf P$_2(\ell,N)$}. The number of isomorphism classes of composition factors of $A[\ell]$ of conductor 1 and $\F_\ell$-dimension at least 2  is  bounded for $A$  in AB$(N)$.

\vspace{5 pt}

\setlength{\hangindent}{0 pt}
\noindent Prop.\! \ref{Onlyell} asserts that {\bf P}$_2(\ell,N)$ necessarily holds for small $\ell$.

\begin{proposition} \label{UseGenBd}
Assume  {\bf P$_1(\ell,N)$} and {\bf P$_2(\ell,N)$}. Then AB$(N)$ is finite.  
\end{proposition}

\proof  We shall denote by $\beta_i$ constants depending only on $\ell$ and $N.$  
Suppose $A$ in AB$(N)$ is good at $\ell$.  Write $\epsilon_0(A[\ell])$ for the number of 1-dimensional constituents of $A[\ell],$  $\gS_\ell(A)$ for the multiset of other irreducible constituents $E$ and   $m_E = \dim_{\F_\ell} E.$ Let $\gS_\ell^1(A)$ be the multiset of those $E$ in $\gS_\ell(A)$ with conductor $N_E = 1.$ We assume that $m_E\le \beta_1$ and $|\gS_\ell^1(A)|\le \beta_2$ as $A$ varies in AB$(N)$.  We have
\begin{equation} \label{GenBd}
2 \dim A = \dim A[\ell] = \epsilon_0(A[\ell]) \; + \sum_{E \in \gS_\ell^1(A)} \hspace{-6 pt} m_E \; + \sum_{N_E > 1}  m_E.
\end{equation}
We bound the last two sums by $\beta_3= \beta_1 (\beta_2+  \Omega(N_A))$, where $\Omega(n)=\sum_{p} \ord_p(n) $.  Theorem 5.3 of \cite{BK2} implies that
$$
\epsilon_0(A[\ell]) \le   2 \, \Omega(N_A) \; +\sum_{E \in \gS_\ell(A)} \hspace{-4 pt} \delta_A(E),
$$
where $\delta_A(E)$ is bounded in terms of $N,$  $\ell$ and the strict class number $h_E$ of $F=\Q(E)$, thanks to \cite[4.3.8, 4.3.13, 4.4.1]{BK2}.
If $m = \dim E$, then 
$
n=[F:\Q]  \le \vert \GL_m(\F_\ell) \vert < \ell^{m^2}.
$ 
Hence the discriminant of $F$ satisfies 
$
\log \vert d_{F/\Q} \vert \le n (2 \log \ell + \log N_E)
$ 
by (\ref{GFdisc}).  An upper bound on the residue of the zeta function \cite{Lou} and a lower bound on the regulator \cite{Sko} show that 
$
h_E \le  \beta_4.$ Thus $\delta_A(E)\le \beta_5$ and 
$\epsilon_0(A[\ell]) \le 2 \Omega(N_A) + \beta_5(\beta_2 + \Omega(N_A)).$   \qed

\medskip

\begin{proposition}[GRH] 
If $N \le 15 683$ is odd and squarefree, then $\dim V$ is bounded for $V$ in ${\rm Irr}(N,2)$ and {\rm AB}$(N)$ is finite.
\end{proposition}

\proof 
Prop.\! \ref{exc2}  restricts the group $G=\Gal(\Q(V)/\Q)$ and Cor. \ref{pdiscSL}, \ref{pdisc} estimate the discriminant of a stem field $K$.  If $\dim V$ is bounded, then the finiteness of ${\rm AB}(N)$ ensues from Prop.\! \ref{UseGenBd}.
Otherwise the root discriminants $\varrho_K$ have an  asymptotic  upper  bound of  $4N^{1/4}$ and  a lower bound of  $8\pi e^{\gamma}$ and so $N > 15 683$.
 \qed

\begin{Rem}
Unconditionally, if $N$ is squarefree,  $V$ is in ${\rm Irr}(N,2)$   and 
$$
2\dim V+1 \ge \max(7,3.06\log N),
$$ 
then $G=\Gal(\Q(V)/\Q)$ can only be $Q$ or $Q\wr \cS_2$ with $Q=O^{\pm}(V)$, $\SP(V)$ or $\SL(V)$.  Otherwise, $G=\cS_m$ or $\cS_m \wr \cS_2$, with $m = 2\dim V+1$ or $2 \dim V+ 2$, as in Rem.\! \ref{SnRep}. By Cor. \ref{pdisc}i, the root discriminant  of a stem field for $\Q(V)$ is at most  $4N^{1/m}$, while it is at least $5.548$ when $m\ge 7$ by \cite{dyd}. 
\end{Rem}

\end{document}